\documentclass[11pt]{article}

\usepackage[margin=1in]{geometry}
\usepackage{amssymb}
\usepackage{amsmath}
\usepackage{amsfonts}
\usepackage{amsthm}
\usepackage{mathtools}
\usepackage{verbatim}
\usepackage[hidelinks]{hyperref}

\numberwithin{equation}{section}
\theoremstyle{plain}
\newtheorem{theorem}{Theorem}[section]
\newtheorem{lemma}[theorem]{Lemma}
\newtheorem{proposition}[theorem]{Proposition}
\newtheorem{corollary}[theorem]{Corollary}
\newtheorem{conjecture}{Conjecture}
\theoremstyle{definition}
\newtheorem{definition}[theorem]{Definition}

\title{On a positivity conjecture of Fasi and Higham for Pad\'e approximants to the matrix exponential}
\author{Francesco Scardino\thanks{Department of Physics, School of Mathematics and Natural Sciences, Bergische Universit\"at Wuppertal, Gau\ss stra\ss e 20, 42119 Wuppertal, Germany
(\texttt{scardino@uni-wuppertal.de}).}}
\date{}

\begin{document}
	
		\maketitle
	
\begin{abstract}
	We study the positivity conjecture of Fasi and Higham \cite[Conjecture~3.2]{higham}, which concerns the Taylor coefficients of the reciprocal of the denominator $q_{k,m}(z)$ of the Pad\'e approximant to the matrix exponential.  The conjecture asserts that these coefficients are positive whenever $m$ is odd.  We prove full positivity for $k=0$ and eventual positivity whenever the unique positive real zero of $q_{k,m}$ is the unique zero of minimal modulus.  We also give an explicit computable sufficient threshold beyond which all coefficients are positive.  Using Hermite scaling with leading-coefficient normalization, we prove that this spectral minimality holds for each fixed odd $m$ and all sufficiently large $k$.
	
	The unrestricted conjecture is nevertheless false.  We give exact rational counterexamples and prove a transient scaling theorem showing that on the scale $n\asymp\sqrt{k+1}$ the normalized coefficients converge to an explicit oscillatory Hermite profile.  For $m=5$ this profile is negative at one point, yielding infinitely many counterexamples.  The counterexamples are not diagonal, and the diagonal case $k=m$, which is the case relevant in practice, remains open.  These results do not imply a failure of the Fasi--Higham algorithm.  Rather, they delimit the range in which the auxiliary positivity conjecture can be used as a theorem.
\end{abstract}
	\medskip
\noindent\textbf{Keywords.} matrix exponential, Pad\'e approximation, reciprocal power series, Hermite polynomials, zero asymptotics

\smallskip
\noindent\textbf{MSC codes.} 41A21, 30B10, 33C45, 65F60, 26C10

\section{Introduction}
Let $f(z)$ be a complex function analytic at $z=0$, let $k\ge0$ be an integer, and let $m\ge1$ be an integer.  The rational function $r_{k,m}(z) = \frac{p_{k,m}(z)}{q_{k,m}(z)}$ is the so-called $[k/m]$ Pad\'e approximant of $f(z)$ at $z=0$ if $p_{k,m}(z)$ and $q_{k,m}(z)$ are polynomials of degree at most $k$ and $m$, respectively, the denominator is normalized so that $q_{k,m}(0) = 1$, and $f(z) - r_{k,m}(z) = \mathcal{O}(z^{k+m+1})$. 
In particular, the numerator and denominator of the $[k/m]$ Pad\'e approximant to the exponential function $f(z)=e^z$ at $0$ are \cite{gautschi2011numerical}
\begin{align}
	&p_{k,m}(z) = \sum_{j=0}^{k} \binom{k}{j} \frac{(k+m-j)!}{(k+m)!} z^j \notag\\
	&q_{k,m}(z) = \sum_{j=0}^{m} (-1)^j \binom{m}{j} \frac{(k+m-j)!}{(k+m)!} z^j\,.\label{eq:qdef}
\end{align}
For later convenience we define the positive coefficients of $q_{k,m}(z)$ by
\begin{equation}\label{eq:omega}
\omega_{k,m}(j)=
\begin{cases}
\displaystyle \binom{m}{j}\frac{(k+m-j)!}{(k+m)!}, & 0\le j\le m,\\
0, & j>m,
\end{cases}
	\qquad j\in\mathbb Z_{\ge0}.
\end{equation}
We use the shorthand
\begin{equation*}
\kappa_{k,m}=\omega_{k,m}(1)=\frac{m}{m+k}\in(0,1].
\end{equation*}
We also write
\begin{equation}\label{eq:cn-def}
	\frac{1}{q_{k,m}(z)} = \sum_{n\ge0} c_n(k,m)\, z^n
\end{equation}
for the Taylor expansion of $q_{k,m}^{-1}(z)$ at the origin, which exists since $q_{k,m}(0)=1$.

Fasi and Higham \cite{higham} developed an algorithm based on Pad\'e approximation for computing the matrix exponential $e^A$, $A\in\mathbb{C}^{d\times d}$, in arbitrary precision.  They also derived a new representation for the truncation error of a general $[k/m]$ Pad\'e approximant and showed how it can be used to compute practical error bounds for truncated Taylor series and diagonal Pad\'e approximants.  In the first case, the bound is cheap to compute, requiring only $O(d^2)$ flops \cite{higham}.  For diagonal Pad\'e approximants the new bound requires $O(d^3)$ low-precision flops, but if their positivity conjecture were true then this cost would be reduced to $O(d^2)$ flops.

We now restate that conjecture.
\begin{conjecture}\label{conj}
Let $k\ge0$ and $m\ge1$ be integers.  If $m$ is odd, then all the coefficients of the Taylor expansion of $q_{k,m}^{-1}(z)$ at $0$ are positive.
\end{conjecture}

The present paper concerns the conjecture as a mathematical statement in its full stated range of parameters.  Our counterexamples do not show that the algorithm of Fasi and Higham fails.  The case $m=1$, which occurs in one of the algorithmic settings, is positive for every $k$, and the explicit finite counterexamples below have $k$ substantially larger than $m$.  Moreover, the $O(d^2)$ bound that the conjecture would provide concerns the diagonal approximants $k=m$, whose denominators are normalized reverse Bessel polynomials \cite{Grosswald1978}.  The diagonal case is precisely the case that our counterexamples leave untouched, and it remains open.  The practical availability of that bound is therefore neither established nor excluded by the present results.  Thus the point of the counterexamples is that the unrestricted positivity conjecture cannot be used without additional hypotheses.  The positive results below identify such hypotheses, while the transient scaling theorem explains the mechanism by which the unrestricted statement fails.

\section{Main results}
The paper has three main parts.  First, we isolate two regimes in which positivity is genuine.  For $k=0$ and every odd $m$ we prove full positivity (Theorem~\ref{thm:main-k0}),
\begin{equation*}
q^{-1}_{0,m}(z)=\sum_{n\ge0}c_n(0,m)z^n,
\qquad c_n(0,m)>0\quad(n\ge0),
\end{equation*}
and for every fixed odd $m$ and every $k$ for which the positive real zero of $q_{k,m}$ is the unique zero of smallest modulus, we prove eventual positivity (Theorem~\ref{thm:eventual-positivity}).  More precisely, if Property~\textnormal{(D1)} in Definition~\ref{def:D1} holds, then there is a computable sufficient index $N_0(k,m)$ such that
\begin{equation*}
c_n(k,m)>0\qquad(n\ge N_0(k,m)).
\end{equation*}
Proposition~\ref{prop:positivity-threshold} gives such an index in terms of the residues and the spectral gap of the poles, and Corollary~\ref{cor:asymptotic-threshold} gives an asymptotic sufficient upper bound for large $k$.

Second, using the Laguerre representation of $q_{k,m}$ and a Hermite scaling limit, we prove that for each fixed odd $m\ge3$ the spectral condition in Property~\textnormal{(D1)} holds for all sufficiently large $k$ (Theorem~\ref{thm:D1-asymptotic-rigorous}).  The proof uses leading-coefficient normalization, rather than normalization at the origin.  This is essential because the limiting Hermite polynomial satisfies $H_m(0)=0$ for odd $m$.  The same asymptotics show that the ratio between the smallest modulus $\rho_{\mathrm{oth}}(k,m)$ among the remaining poles and the positive pole modulus $R_{k,m}$ is
\begin{equation*}
\frac{\rho_{\mathrm{oth}}(k,m)}{R_{k,m}}=1+\frac{\delta_m}{3(k+1)}+O((k+1)^{-3/2}),
\end{equation*}
where $H_m$ is the physicists' Hermite polynomial and $\delta_m$ is its smallest nonzero squared zero.  Thus the pole-dominance argument provides a sufficient threshold on the scale $n=O(k)$.

Third, we show that the unrestricted conjecture is nevertheless false.  We give exact rational finite counterexamples, verified by a short reproducible symbolic computation, namely $c_{38}(66,5)<0$, $c_{37}(44,7)<0$, and $c_{42}(45,9)<0$ (Theorem~\ref{thm:counterexamples-again}).  We also prove a transient scaling theorem for $n\asymp\sqrt{k+1}$ (Theorem~\ref{thm:transient-scaling}).  With $N=k+1$ and $n_N/\sqrt N\to\lambda$, the normalized coefficients converge to an explicit oscillatory Hermite profile $\Psi_m(\lambda)$, and for $m=5$ this profile is negative at $\lambda=5$.  Consequently there are infinitely many pairs $(k,n)$ for which $c_n(k,5)<0$ (Corollary~\ref{cor:m5-infinite-counterexamples}).  This negative window on the scale $n\asymp\sqrt{k}$ is compatible with the sufficient eventual-positivity threshold on the larger scale $n=O(k)$.

None of these results settles the diagonal case $k=m$, which is the case relevant to the error bounds of \cite{higham}.  The counterexamples are not diagonal, and Theorem~\ref{thm:D1-asymptotic-rigorous} concerns fixed $m$ and large $k$.  Diagonal positivity therefore remains an open problem.

The paper is organized as follows.  Section~\ref{sec:firststep} proves the case $k=0$.  Section~\ref{sec:structural} introduces the general odd-$m$ notation and the representation formulas and proves that all zeros are simple.  Section~\ref{app:coeff} gives the Enestr\"om--Kakeya bounds and the elementary zero-location and monotonicity facts, including the unique positive real zero for odd $m$.  Sections~\ref{sec:asymptotics} and \ref{sec:rigorous_proof} prove the Hermite zero asymptotics and Property~\textnormal{(D1)} for large $k$.  Sections~\ref{sec:failure} and \ref{sec:transient} use partial fractions to obtain the positivity threshold, the exact counterexamples, and the transient scaling limit.  Section~\ref{sec:conclusion} summarizes the results.

\section{The case $k=0$}\label{sec:firststep}

The first positive result settles the case $k=0$ completely.
\begin{theorem}\label{thm:main-k0}
Let $m\ge1$ be an odd integer.  Then $c_n(0,m)>0$ for every $n\ge0$.
\end{theorem}

The polynomial $q_{k,m}(z)$ is simplest when $k=0$.  Indeed, by the second line of \eqref{eq:qdef},
\begin{equation*}
		q_{0,m}(z)= \sum_{j=0}^{m} \frac{(-1)^j}{j!} z^j\,.
\end{equation*}
We write $q_{0,m}(z)$ in terms of the upper incomplete gamma function
\begin{equation*}
	\Gamma(s,z) = \int_z^\infty t^{s-1} e^{-t}\,dt\,,
\end{equation*}
which, for an integer $m\ge1$, satisfies \cite{NIST:DLMF}
\begin{equation*}	\Gamma(m,z) = (m-1)!e^{-z}\sum_{j=0}^{m-1}\frac{z^j}{j!}\,,
\end{equation*}
and so
\begin{equation}
	\label{eq:resum}
	q_{0,m}(z) = \frac{e^{-z}}{m!}\,\Gamma (m+1,-z)\,.
\end{equation}
We begin by showing that the Taylor expansion of $q^{-1}_{0,m}(z)$ has positive coefficients when $m$ is odd.  Its radius of convergence is at least $1$, as shown in Section~\ref{app:coeff}.
To expand the reciprocal in \eqref{eq:resum}, we use the series \cite{NIST:DLMF}
\begin{equation*}
	\Gamma(s,x)=\Gamma(s)-\sum_{n=0}^{\infty}\frac{(-1)^n}{n!}\frac{x^{s+n}}{s+n}\,,
\end{equation*}
where $\Gamma(s)$ is the ordinary gamma function.  Setting $s=m+1$ and $x=-z$ gives
\begin{equation*}
	\Gamma(m+1,-z)=m!+(-1)^m\sum_{n=0}^{\infty}\frac{1}{n!}\frac{z^{n+m+1}}{n+m+1}\,.
\end{equation*}
Therefore, $q^{-1}_{0,m}(z)$ becomes
\begin{equation*}
	q^{-1}_{0,m}(z) = \frac{e^{z}  }{1+\frac{(-1)^m}{m!}\sum_{n=0}^{\infty}\frac{z^{n+m+1}}{n!(n+m+1)}}\,.
\end{equation*}
The series in the denominator has positive coefficients and we call it
\begin{equation*}
	F_{m}(z) = \frac{1}{m!}\sum_{n=0}^{\infty}\frac{z^{n+m+1}}{n!(n+m+1)}\,,
\end{equation*} 
and so
\begin{equation}
	\label{eq:key}
	q^{-1}_{0,m}(z) = \frac{e^{z}  }{1+(-1)^mF_m(z)}\,.
\end{equation}
If we set $m=2\ell+1$ odd, the denominator gets a minus sign.  For $|F_{2\ell+1}(z)|<1$ the geometric series converges, and the identity holds coefficientwise as a formal Taylor series at the origin:
\begin{equation*}
	q^{-1}_{0,2\ell+1}(z) = \frac{e^{z}}{1-F_{2\ell+1}(z)}=e^{z}\sum_{n=0}^\infty (F_{2\ell+1}(z))^n\,.
\end{equation*}
This product has positive Taylor coefficients: every $(F_{2\ell+1}(z))^n$ has nonnegative coefficients, and the $n=0$ term contributes $e^z$, whose coefficients are all positive.  It is also immediate from \eqref{eq:key} why the same argument does not apply for even $m$.
We conclude that, for every $\ell\ge0$,
\begin{equation*}
		q^{-1}_{0,2\ell+1}(z) = \sum_{n=0}^\infty c_n(0,2\ell+1)\, z^n,\qquad c_n(0,2\ell+1) >0\quad(n\ge0),
\end{equation*}
where the positivity follows from the displayed factorization.  This proves Theorem~\ref{thm:main-k0}.

\section{Representations and notation for general odd $m$}\label{sec:structural}\label{app:polynomials}

Fix an \emph{odd} integer $m\ge3$ and an integer $k\ge0$.  In the notation of \eqref{eq:qdef} and \eqref{eq:cn-def}, the coefficients \eqref{eq:omega} admit the product form
\begin{equation*}
\omega_{k,m}(j)=\frac{1}{j!}\prod_{t=0}^{j-1}\frac{m-t}{m+k-t},
\end{equation*}
which will be used repeatedly below.  We discuss only odd values $m\ge3$ because for $m=1$ the conjecture is immediate from
\begin{equation*}
q_{k,1}(z)=1-\frac{1}{1+k}z
\end{equation*}
and $q^{-1}_{k,1}(z)$ is a geometric series.

The terminating hypergeometric representation is \cite{SzegoOP}
\begin{equation*}
		q_{k,m}(z) = \sum_{j=0}^{m} (-1)^j \binom{m}{j} \frac{(k+m-j)!}{(k+m)!} z^j
	= {}_1F_1\!\left(\left.
	\begin{matrix}-m\\-m-k\end{matrix}\right|-z\right)\,.
\end{equation*}
Using the factorial integral
\begin{equation*}
	(k+m-j)! = \int_{0}^{\infty} t^{m+k-j}\, e^{-t}\, dt\,,
\end{equation*}
we obtain
\begin{align*}
	q_{k,m}(z) = \int_{0}^{\infty} \frac{t^{m+k}}{(m+k)!}e^{-t}\sum_{j=0}^m(-1)^j\binom{m}{j}(t^{-1}z)^j\,dt,
\end{align*}
which becomes
\begin{align}\label{eq:integralrep}
	q_{k,m}(z) = \frac{1}{(m+k)!} \int_{0}^{\infty}t^k\,(t-z)^m e^{-t}\,dt\,.
\end{align}
The identity \eqref{eq:integralrep} is valid for every integer $m\ge1$.  Its derivation does not use a parity assumption.
For $\operatorname{Re}z>0$, expanding $(t-1)^m$ and evaluating the resulting gamma integrals also gives
\begin{align*}
	q_{k,m}(z) = \frac{z^{k+m+1}}{(m+k)!} \int_{0}^{\infty}t^{k}\,(t-1)^m\,e^{-z\,t}\,dt\,.
\end{align*}
There is a known relation \cite{SzegoOP} between the confluent hypergeometric function and the generalized Laguerre polynomials, namely
\begin{equation*}
	L^{(\alpha)}_n(x) =\frac{(\alpha+1)_n} {n!}
	{}_1F_1\!\left(\left.
	\begin{matrix}-n\\\alpha+1\end{matrix}\right|x\right),
\end{equation*}
where $(\alpha+1)_n$ is the Pochhammer symbol. Then we can write
\begin{equation}
	\label{eq:laguerre_rep}
q_{k,m}(z)
	= \frac{m!}{(-m-k)_m}L^{(-m-k-1)}_m(-z).
\end{equation}

\begin{lemma}\label{lem:simple-zeros}
For every $m\ge1$ and $k\ge0$, all zeros of $q_{k,m}$ are simple.
\end{lemma}

\begin{proof}
By \eqref{eq:laguerre_rep}, it suffices to prove simplicity for
$y(x)=L_m^{(-m-k-1)}(x)$.  This polynomial satisfies
\begin{equation*}
xy''(x)+(-m-k-x)y'(x)+m y(x)=0.
\end{equation*}
If $y(x_0)=y'(x_0)=0$ at an ordinary point $x_0\ne0$, uniqueness for
linear differential equations forces $y\equiv0$, a contradiction.  Finally,
\begin{equation*}
L_m^{(-m-k-1)}(0)=\binom{-k-1}{m}
=(-1)^m\binom{k+m}{m}\ne0,
\end{equation*}
so $0$ is not a zero.
\end{proof}

\section{Basic zero location and real-root structure}\label{app:coeff}

We begin this section by restating a generalization of the Enestr\"om--Kakeya theorem
\cite{Gardner2014,melman2018,Enest1,Enest2,Kak}.

\begin{theorem}
	All the zeros of the real polynomial $p(z) = \sum_{j=0}^{n} a_j z^j $ with positive coefficients lie in the annulus
	\begin{equation*}
	\left\{z \in \mathbb{C} :
	\min_{0 \leq j \leq n-1} \frac{a_j}{a_{j+1}}
	\leq |z| \leq
	\max_{0 \leq j \leq n-1} \frac{a_j}{a_{j+1}}\right\}.
	\end{equation*}
\end{theorem}

This applies precisely to the case
\begin{equation*}
q_{k,m}(-z)=\sum_{j=0}^{m} \omega_{k,m}(j) z^j,
\end{equation*}
since the zeros of a polynomial $p(z)$ are the zeros of $p(-z)$ reflected through the origin.
We obtain directly from the definition of the coefficients in \eqref{eq:omega} that
\begin{equation*}	\frac{\omega_{k,m}(j)}{\omega_{k,m}(j+1)}=(j+1)\frac{m-j+k}{m-j}\ge1,
\end{equation*}
with equality only for $k=0$ and $j=0$, so $\omega_{k,m}(0)\ge\omega_{k,m}(1)>\cdots>\omega_{k,m}(m)>0$.
Therefore the Enestr\"om--Kakeya theorem gives
\begin{equation*}	\frac{m+k}{m}\leq |z|\leq m(k+1)\,.
\end{equation*}
Here the standard theorem gives a nonstrict lower inequality, while part (b)
of Theorem~3.2 in \cite{melman2018}, with
$\delta_1=m/(m+k)$ and $\delta_0=0$, makes it strict in the present range
$m\ge3$.  In particular, all zeros of $q_{k,m}$ lie outside the unit disk,
and $q_{k,m}^{-1}$ is holomorphic on a neighborhood of the closed disk
$|z|\leq(m+k)/m$.
\subsection{Derivative identity, strict decrease on $\mathbb{R}$, and uniqueness of the positive zero}

\begin{lemma}\label{lem:derivative}
	For all integers $m\ge1$ and $k\ge0$,
	\begin{equation*}
	q'_{k,m}(z)\ =\ -\,\kappa_{k,m}\,q_{k,m-1}(z).
	\end{equation*}
\end{lemma}

\begin{proof}
	Compute directly from the coefficients:
	\begin{equation*}
	q'_{k,m}(z)=\sum_{j=1}^{m}(-1)^j j\,\omega_{k,m}(j)\,z^{j-1}.
	\end{equation*}
	Since
	\begin{align*}
		j\,\omega_{k,m}(j)
		&=
		j\cdot\frac{1}{j!}\prod_{t=0}^{j-1}\frac{m-t}{m+k-t}
		=
		\frac{m}{m+k}\cdot \frac{1}{(j-1)!}\prod_{t=1}^{j-1}\frac{m-t}{m+k-t}\\
		&=\kappa_{k,m}\,\omega_{k,m-1}(j-1),
	\end{align*}
	we get
	\begin{equation*}
	q'_{k,m}(z)
	=
	-\kappa_{k,m}\,\sum_{j=1}^m(-1)^{j-1}\omega_{k,m-1}(j-1)z^{j-1}
	=
	-\kappa_{k,m}\,q_{k,m-1}(z).
	\end{equation*}
\end{proof}

\begin{lemma}\label{lem:monotone}
	Let $m\ge3$ be odd and $k\ge0$. Then $q_{k,m}$ is strictly decreasing on $\mathbb{R}$. In particular, $q_{k,m}$ has a unique real zero $R_{k,m}>0$, and $q'_{k,m}(R_{k,m})<0$.
\end{lemma}

\begin{proof}
	By Lemma~\ref{lem:derivative}, $q'_{k,m}(x)=-\kappa_{k,m}\,q_{k,m-1}(x)$.
	For odd $m$, the integer $m-1$ is even.  The integral representation from Section~\ref{app:polynomials} gives
	\begin{equation*}
	q_{k,m-1}(x)=\sum_{j=0}^{m-1}(-1)^j\omega_{k,m-1}(j)x^j
	\end{equation*}
	and this polynomial admits the representation
	\begin{equation*}
	q_{k,m-1}(x)=\frac{1}{(m+k-1)!}\int_0^\infty t^k\, (t-x)^{m-1}\,e^{-t}\,dt.
	\end{equation*}
	Here $(t-x)^{m-1}\ge0$ for all $t,x\in\mathbb{R}$ because $m-1$ is even, and the weight $e^{-t}t^k$ is positive on $(0,\infty)$, so $q_{k,m-1}(x)>0$ for all $x\in\mathbb{R}$.
	Thus $q'_{k,m}(x)<0$ for all $x\in\mathbb{R}$.  Since $q_{k,m}(0)=1$, while its odd degree and negative leading coefficient give $q_{k,m}(x)\to-\infty$ as $x\to+\infty$ and $q_{k,m}(x)\to+\infty$ as $x\to-\infty$, it has a unique real zero $R_{k,m}>0$, and $q'_{k,m}(R_{k,m})<0$.
\end{proof}
\section{Asymptotic clustering and spectral minimality}\label{sec:asymptotics}
For a polynomial $p$, let $\mathcal{Z}(p)$ denote its set of zeros.  We begin with the definition used throughout the rest of the paper.
\begin{definition}\label{def:D1}
	Let $m\ge3$ be odd and $k\ge0$. We say that the pair $(k,m)$ satisfies \emph{Property \textnormal{(D1)}} if every zero $w\in\mathcal{Z}(q_{k,m})$ other than the positive real zero $R_{k,m}$ satisfies
	\begin{equation*}
	|w|>R_{k,m}.
	\end{equation*}
	Equivalently, $R_{k,m}$ is the unique zero of minimal modulus.
\end{definition} 
In this section we describe the large-$k$ behaviour of the zeros of $q_{k,m}$ for
fixed odd $m\ge3$ using only the explicit coefficients of $q_{k,m}$. We first
pass to a scaling limit showing that all zeros cluster near $k$ and that the spectral gap
between the dominant real zero and the remaining zeros closes as $k\to\infty$.
The refined Hermite asymptotics of Section~\ref{sec:rigorous_proof} will then identify the unique zero of minimal modulus and quantify the spectral gap between the positive real zero and the complex zeros.  This provides the zero-scale context for both the eventual-positivity result and the transient counterexamples.

Throughout this section, for given $m\ge1$ and $k\ge0$ we denote the zeros of
$q_{k,m}$ (counted with multiplicity) by
\begin{equation*}
\zeta_1(k,m),\dots,\zeta_m(k,m).
\end{equation*}
When $m$ is odd, Lemma~\ref{lem:monotone} gives a unique real zero, and we
label it $R_{k,m}=\zeta_1(k,m)>0$.

We now pass to a scaling limit at the $k$--scale.

\begin{theorem}\label{thm:scaling-limit}
	Fix $m\ge1$. For each $k\ge1$ define
	\begin{equation*}
	P_k(x) = q_{k,m}(k x).
	\end{equation*}
	Then, as $k\to\infty$,
	\begin{equation*}
	P_k(x) \longrightarrow (1-x)^m
	\end{equation*}
	locally uniformly on $\mathbb{C}$. In particular, for every compact set
	$S \subset \mathbb{C}$,
	\begin{equation*}
	\sup_{x\in S} \bigl| P_k(x) - (1-x)^m \bigr| \longrightarrow 0
	\qquad (k\to\infty).
	\end{equation*}
\end{theorem}

\begin{proof}
	From \eqref{eq:qdef} we have
	\begin{equation*}
	q_{k,m}(z)
	= \sum_{j=0}^{m} (-1)^j \binom{m}{j} \frac{(k+m-j)!}{(k+m)!} z^j.
	\end{equation*}
	Substituting $z = k x$ gives
	\begin{equation}\label{eq:Pk-expansion}
		P_k(x)
		= \sum_{j=0}^{m} (-1)^j \binom{m}{j} a_{k,j}\, x^j,
	\end{equation}
	where
	\begin{equation*}
	a_{k,j}
	= \frac{(k+m-j)!}{(k+m)!}\, k^j
	= \frac{k^j}{(k+m)(k+m-1)\cdots (k+m-j+1)}.
	\end{equation*}
	For $j=0$ we have $a_{k,0}=1$. For $1\le j\le m$,
	\begin{equation*}
	a_{k,j}
	= \prod_{r=0}^{j-1} \frac{k}{k+m-r}
	= \prod_{r=0}^{j-1} \frac{1}{1 + \frac{m-r}{k}}.
	\end{equation*}
	Fix $j$ and let $k\to\infty$. Each factor tends to $1$, and the convergence
	is uniform in $j$ for $0\le j\le m$, since
	\begin{equation*}
	\left| \frac{1}{1 + \frac{m-r}{k}} - 1 \right|
	= \frac{m-r}{k + m - r}
	\le \frac{m}{k}
	\end{equation*}
	for all $0\le r\le j-1 \le m-1$. Hence there exists $C_m>0$ such that
	\begin{equation*}
	|a_{k,j} - 1|
	\le \Bigl(1 + \frac{m}{k}\Bigr)^m - 1
	\le C_m \frac{1}{k}
	\end{equation*}
	for all $0\le j\le m$, and $a_{k,j} \to 1$ uniformly in $j$ as $k\to\infty$.
	
	Let $S\subset\mathbb{C}$ be compact and set
	\begin{equation*}
	R = \max\{1, \sup_{x\in S}|x|\}.
	\end{equation*}
	Using \eqref{eq:Pk-expansion} we obtain
	\begin{equation*}
	\sup_{x\in S} |P_k(x) - (1-x)^m|
	\le \sum_{j=0}^m \binom{m}{j} |a_{k,j} - 1| R^j
	\le C_m' \max_{0\le j\le m} |a_{k,j} - 1|,
	\end{equation*}
	for some constant $C_m'>0$ depending only on $m$ and $R$. Since
	$\max_j|a_{k,j}-1|\to0$ as $k\to\infty$, the right-hand side tends to $0$.
	This proves local uniform convergence on $\mathbb{C}$.
\end{proof}

The scaling limit immediately yields clustering of the zeros.

\begin{corollary}\label{cor:clustering}
	Fix $m\ge1$ and let $k\ge1$. Let $\zeta_1(k,m),\dots,\zeta_m(k,m)$ be the zeros of $q_{k,m}$ and
	set
	\begin{equation*}
	x_j(k,m) = \frac{\zeta_j(k,m)}{k},\qquad 1\le j\le m.
	\end{equation*}
	Then for every $\varepsilon>0$ there exists $K(\varepsilon,m)$ such that
	for all $k\ge K(\varepsilon,m)$,
	\begin{equation*}
	|x_j(k,m) - 1| < \varepsilon
	\qquad\text{for all } 1\le j\le m.
	\end{equation*}
	In particular, for all sufficiently large $k$,
	\begin{equation*}
	(1-\varepsilon) k < |\zeta_j(k,m)| < (1+\varepsilon) k
	\qquad\text{for all } 1\le j\le m.
	\end{equation*}
\end{corollary}

\begin{proof}
	By Theorem~\ref{thm:scaling-limit}, $P_k(x)=q_{k,m}(k x)$ converges locally
	uniformly to $(1-x)^m$. The latter polynomial has a single zero at $x=1$ of
	multiplicity $m$. Let $\varepsilon>0$ and consider the disk
	\begin{equation*}
	D_\varepsilon = \{ x\in\mathbb{C} : |x-1| < \varepsilon \}.
	\end{equation*}
	On the boundary $\partial D_\varepsilon$ we have $(1-x)^m\neq0$, and the
	local uniform convergence implies that for $k$ large enough
	$P_k(x)$ is uniformly close to $(1-x)^m$ on $\partial D_\varepsilon$.
	By Rouch\'e's theorem, for all sufficiently large $k$ the polynomials
	$P_k$ and $(1-x)^m$ have the same number of zeros (counted with multiplicity)
	inside $D_\varepsilon$, namely $m$. Since $\deg P_k = m$, this means that
	all zeros $x_j(k,m)$ of $P_k$ lie in $D_\varepsilon$ for $k$ large.
	
	The stated inequalities for $|\zeta_j(k,m)| = k |x_j(k,m)|$ follow immediately.
\end{proof}

As an immediate consequence, both the positive zero and the closest remaining
zero have modulus asymptotic to $k$.

\begin{corollary}\label{cor:linear-growth}
	Let $m\ge3$ be odd. Let $R_{k,m}>0$ be the unique positive real zero of $q_{k,m}$
	from Lemma~\ref{lem:monotone} and set
	\begin{equation*}
	\rho_{\mathrm{oth}}(k,m) = \min_{2\le j\le m} |\zeta_j(k,m)|
	\end{equation*}
	for the smallest modulus among the remaining zeros. Then, as $k\to\infty$,
	\begin{equation*}
	\frac{R_{k,m}}{k} \longrightarrow 1,
	\qquad
	\frac{\rho_{\mathrm{oth}}(k,m)}{k} \longrightarrow 1,
	\end{equation*}
	and in particular
	\begin{equation*}
	\lim_{k\to\infty} \frac{\rho_{\mathrm{oth}}(k,m)}{R_{k,m}} = 1.
	\end{equation*}
\end{corollary}

\begin{proof}
	For each $k$, the set $\{\zeta_j(k,m)\}_{j=1}^m$ is exactly the zero set of
	$q_{k,m}$. By Corollary~\ref{cor:clustering}, for every $\varepsilon>0$ and
	all sufficiently large $k$ we have, for all $j$,
	\begin{equation*}
	(1-\varepsilon)k < |\zeta_j(k,m)| < (1+\varepsilon)k.
	\end{equation*}
	Applying this bound to $j=1$ gives the corresponding inequalities for
	$R_{k,m}$, and taking the minimum over $j\ge2$ gives them for
	$\rho_{\mathrm{oth}}(k,m)$. Fixing $\varepsilon>0$, dividing by $k$, and letting $k\to\infty$
	yields
	\begin{equation*}
	1-\varepsilon \le \liminf_{k\to\infty} \frac{R_{k,m}}{k}
	\le \limsup_{k\to\infty} \frac{R_{k,m}}{k} \le 1+\varepsilon,
	\end{equation*}
	and likewise for $\rho_{\mathrm{oth}}(k,m)/k$. Since $\varepsilon>0$ is arbitrary, both
	limits equal $1$, and therefore
	\begin{equation*}
	\frac{\rho_{\mathrm{oth}}(k,m)}{R_{k,m}}
	= \frac{\rho_{\mathrm{oth}}(k,m)/k}{R_{k,m}/k}
	\longrightarrow 1.
	\end{equation*}
\end{proof}
	Using the Enestr\"om--Kakeya bounds from Section~\ref{app:coeff}, one can
	already obtain cruder estimates of the form
	\begin{equation*}
	C_{m,-}k \le \min_j |\zeta_j(k,m)| \le \max_j |\zeta_j(k,m)| \le C_{m,+}k
\end{equation*}
	for suitable constants $0<C_{m,-}\le C_{m,+}<\infty$. Corollary~\ref{cor:clustering}
	strengthens these bounds by showing that, after scaling by $k$, \emph{all}
	zeros converge to the single point $1$.
	The refined scaling below determines which branch has the smallest modulus
	and quantifies the closing spectral gap.

\section{Asymptotic proof of Property (D1)}\label{sec:rigorous_proof}

We now prove the refined large-$k$ zero asymptotics used to establish Property~\textnormal{(D1)}.  Throughout this section $m\ge3$ is fixed and odd, and
\begin{equation*}
N=k+1.
\end{equation*}
By \eqref{eq:laguerre_rep}, the zeros of $q_{k,m}$ are the negatives of the zeros of
\begin{equation*}
y(x)=L_m^{(-m-k-1)}(x)=L_m^{(-m-N)}(x).
\end{equation*}
The Laguerre differential equation is
\begin{equation}\label{eq:Laguerre-ODE-N-revised}
xy''(x)+(1-m-N-x)y'(x)+m y(x)=0.
\end{equation}
Set
\begin{equation*}x=-N+\sqrt{2N}\,h,
\qquad \varepsilon=N^{-1/2},
\end{equation*}
and write $Y(h,\varepsilon)=y(-N+\sqrt{2N}h)$.  Substitution in \eqref{eq:Laguerre-ODE-N-revised} gives the exact equation
\begin{equation}\label{eq:Y-ODE-expanded-h-revised}
Y_{hh}+2hY_h-2mY-
\sqrt2\,\varepsilon\bigl(hY_{hh}+(1-m)Y_h\bigr)=0.
\end{equation}

A normalization at the origin, $Y(0,\varepsilon)=1$, would be singular for odd $m$, because the limiting Hermite polynomial vanishes at the origin.  We therefore normalize by the leading coefficient.  Put $z=ih$ and define
\begin{equation*}
\widetilde Y(z,\varepsilon)=c(\varepsilon)Y(-iz,\varepsilon),
\end{equation*}
where the nonzero scalar $c(\varepsilon)$ is chosen so that $[z^m]\widetilde Y(z,\varepsilon)=2^m$.  Since the differential equation is homogeneous, \eqref{eq:Y-ODE-expanded-h-revised} becomes
\begin{equation}\label{eq:Y-ODE-z-eps-revised}
\widetilde Y_{zz}-2z\widetilde Y_z+2m\widetilde Y
+i\sqrt2\,\varepsilon\bigl(z\widetilde Y_{zz}+(1-m)\widetilde Y_z\bigr)=0,
\end{equation}
with the normalization
\begin{equation}\label{eq:leading-normalization}
[z^m]\widetilde Y(z,\varepsilon)=2^m.
\end{equation}
This normalization is nonsingular, since the physicists' Hermite polynomial satisfies $[z^m]H_m(z)=2^m$.

Write
\begin{equation*}
\widetilde Y(z,\varepsilon)=\sum_{r=0}^m a_r(\varepsilon)z^r,
\qquad a_m(\varepsilon)=2^m,
\qquad a_{m+1}(\varepsilon)=a_{m+2}(\varepsilon)=0.
\end{equation*}
Equating coefficients in \eqref{eq:Y-ODE-z-eps-revised} gives, for $s=0,\ldots,m$,
\begin{equation*}(s+2)(s+1)a_{s+2}(\varepsilon)+2(m-s)a_s(\varepsilon)
+i\sqrt2\,\varepsilon(s+1)(s+1-m)a_{s+1}(\varepsilon)=0.
\end{equation*}
This recurrence determines the coefficients downward from $a_m(\varepsilon)=2^m$, and hence each $a_s(\varepsilon)$ is analytic at $\varepsilon=0$.  At $\varepsilon=0$ the recurrence is the Hermite coefficient recurrence with leading coefficient $2^m$, so
\begin{equation}\label{eq:Hermite-limit-revised}
\widetilde Y(z,\varepsilon)=H_m(z)+\varepsilon\widetilde Y_1(z)+O(\varepsilon^2)
\end{equation}
coefficientwise, hence locally uniformly on compact sets.

Let
\begin{equation*}
\mathcal L[f]=f''-2zf'+2mf.
\end{equation*}
Collecting the $\varepsilon$ terms in \eqref{eq:Y-ODE-z-eps-revised} yields
\begin{equation*}\mathcal L[\widetilde Y_1]
=-i\sqrt2\bigl(zH_m''+(1-m)H_m'\bigr).
\end{equation*}
Using the standard Hermite identities
\begin{equation*}
H_m'=2mH_{m-1},
\qquad H_m''=4m(m-1)H_{m-2},
\qquad zH_n=\frac12H_{n+1}+nH_{n-1},
\end{equation*}
the source simplifies to
\begin{equation*}
zH_m''+(1-m)H_m'=4m(m-1)(m-2)H_{m-3}.
\end{equation*}
Since $\mathcal L[H_r]=2(m-r)H_r$ and the normalization \eqref{eq:leading-normalization} forces $[z^m]\widetilde Y_1=0$, there is no homogeneous $H_m$ component and
\begin{equation*}\widetilde Y_1(z)=
-i\sqrt2\,\frac{2}{3}m(m-1)(m-2)H_{m-3}(z).
\end{equation*}

By coefficientwise convergence and the simplicity of the Hermite zeros, Rouch\'e's theorem gives one zero branch of $\widetilde Y$ near each zero of $H_m$.  Let $\xi$ be such a zero.  The corresponding zero branch has the form
\begin{equation*}
z_\xi(\varepsilon)=\xi+\varepsilon\delta_\xi+O(\varepsilon^2).
\end{equation*}
From \eqref{eq:Hermite-limit-revised},
\begin{equation*}
\delta_\xi=-\frac{\widetilde Y_1(\xi)}{H_m'(\xi)}.
\end{equation*}
Using $H_m'(\xi)=2mH_{m-1}(\xi)$ and the recurrence relation evaluated at the zero $\xi$,
\begin{equation*}
H_{m-3}(\xi)=\frac{2\xi^2-m+1}{2(m-1)(m-2)}H_{m-1}(\xi),
\end{equation*}
we get
\begin{equation*}\delta_\xi=\frac{i\sqrt2}{6}(2\xi^2-m+1).
\end{equation*}
Because $z=ih$, the zeros $\ell_{\xi,N}$ of $L_m^{(-m-k-1)}$ satisfy
\begin{equation*}\ell_{\xi,N}
=-N-i\sqrt{2N}\,\xi+\frac{2\xi^2-m+1}{3}+O(N^{-1/2}),
\end{equation*}
uniformly over the finite zero set of $H_m$.  Therefore the zeros $\zeta_{\xi,N}=-\ell_{\xi,N}$ of $q_{k,m}$ satisfy
\begin{equation}\label{eq:q-zero-asymptotic-revised}
\zeta_{\xi,N}
=N+i\sqrt{2N}\,\xi-\frac{2\xi^2-m+1}{3}+O(N^{-1/2}).
\end{equation}
Let $x_{0,N}$ be the scaled zero on the branch corresponding to the Hermite zero $\xi=0$, so that $\zeta_{0,N}=N+i\sqrt{2N}\,x_{0,N}$.  Because $q_{k,m}$ has real coefficients, conjugation sends a scaled zero $x$ to $-\overline{x}$.  The branch near $0$ is unique, hence $x_{0,N}=-\overline{x_{0,N}}$ and $\zeta_{0,N}$ is real.  Equation~\eqref{eq:q-zero-asymptotic-revised} makes this zero positive for large $N$, so Lemma~\ref{lem:monotone} identifies it with $R_{k,m}$.  Thus
\begin{equation}\label{eq:R-asymptotic-revised}
R_{k,m}=N+\frac{m-1}{3}+O(N^{-1/2}).
\end{equation}
For any nonzero Hermite root $\xi$, comparison of moduli gives
\begin{equation}\label{eq:modulus-gap-revised}
|\zeta_{\xi,N}|^2-R_{k,m}^2
=\frac23N\xi^2+O(N^{1/2})>0
\end{equation}
for all sufficiently large $N$.  This proves the desired spectral minimality.

\begin{theorem}\label{thm:D1-asymptotic-rigorous}
Let $m\ge3$ be odd.  Then there exists $K_m$ such that, for all $k\ge K_m$, the unique positive real zero $R_{k,m}$ of $q_{k,m}$ is the unique zero of minimal modulus.  Equivalently, Property~\textnormal{(D1)} holds for all sufficiently large $k$.
\end{theorem}

\begin{proof}
The branch identification preceding \eqref{eq:R-asymptotic-revised} gives the unique positive real zero.  For every other branch $\xi\ne0$, \eqref{eq:modulus-gap-revised} is positive for all sufficiently large $N=k+1$, uniformly over the finite set of nonzero Hermite roots.  Hence $|\zeta_{\xi,N}|>R_{k,m}$ for every other zero once $k\ge K_m$.
\end{proof}

\section{Eventual positivity and failure of full positivity}\label{sec:failure}

We now connect the location of the zeros with the Taylor coefficients of $q_{k,m}^{-1}$.  The first result is the standard partial fraction reduction, stated in a form used both for eventual positivity and for the transient scaling theorem.

\begin{lemma}\label{lem:partial-fractions}
Let $m\ge1$ and $k\ge0$, and let $\zeta_1,\ldots,\zeta_m$ be the zeros of $q_{k,m}$, which are simple by Lemma~\ref{lem:simple-zeros}.  Then
\begin{equation*}\frac{1}{q_{k,m}(z)}=\sum_{j=1}^m\frac{\mathcal R_j}{1-z/\zeta_j},
\qquad
\mathcal R_j=-\frac{1}{\zeta_jq_{k,m}'(\zeta_j)}.
\end{equation*}
Consequently, for the Taylor coefficients of $q_{k,m}^{-1}$,
\begin{equation}\label{eq:cn-sum-again}
c_n(k,m)=\sum_{j=1}^m \mathcal R_j\zeta_j^{-n}\qquad(n\ge0).
\end{equation}
\end{lemma}

\begin{proof}
Write $q_{k,m}(z)=(z-\zeta_j)g_j(z)$ near a zero $\zeta_j$.  Then
\begin{equation*}
\frac1{g_j(z)}=\frac1{g_j(\zeta_j)}+
\left(\frac1{g_j(z)}-\frac1{g_j(\zeta_j)}\right),
\end{equation*}
where the expression in parentheses vanishes at $z=\zeta_j$.  Hence the principal part of $1/q_{k,m}$ at $\zeta_j$ is
\begin{equation*}
\frac{1}{q_{k,m}'(\zeta_j)}\frac1{z-\zeta_j}
=-\frac{1}{\zeta_jq_{k,m}'(\zeta_j)}\frac1{1-z/\zeta_j}.
\end{equation*}
Subtracting the sum of these principal parts gives an entire function that tends to $0$ at infinity, and hence is identically zero by Liouville's theorem.  Expanding each $(1-z/\zeta_j)^{-1}$ in the disk $|z|<\min_j|\zeta_j|$ gives \eqref{eq:cn-sum-again}.
\end{proof}

\begin{theorem}\label{thm:eventual-positivity}
Fix odd $m\ge3$ and $k\ge0$.  Suppose that the pair $(k,m)$ satisfies Property~\textnormal{(D1)}, and let $R_{k,m}$ be the positive real zero of $q_{k,m}$.  Then there exists $N_0=N_0(k,m)$ such that
\begin{equation*}
c_n(k,m)>0\qquad(n\ge N_0).
\end{equation*}
\end{theorem}

\begin{proof}
Label the zeros so that $\zeta_1=R_{k,m}$.  By Lemma~\ref{lem:monotone}, $q_{k,m}'(R_{k,m})<0$, hence
\begin{equation*}
\mathcal R_1=-\frac1{R_{k,m}q_{k,m}'(R_{k,m})}>0.
\end{equation*}
Let $\rho_{\mathrm{oth}}=\min_{j\ge2}|\zeta_j|$.  Property~\textnormal{(D1)} gives $\rho_{\mathrm{oth}}>R_{k,m}$.  From Lemma~\ref{lem:partial-fractions},
\begin{equation*}
c_n(k,m)=\mathcal R_1R_{k,m}^{-n}+\varepsilon_n,
\qquad
|\varepsilon_n|\le \mathcal B\rho_{\mathrm{oth}}^{-n},
\qquad
\mathcal B=\sum_{j=2}^m |\mathcal R_j|.
\end{equation*}
With $\gamma=\rho_{\mathrm{oth}}/R_{k,m}>1$, choose $N_0$ so that $\mathcal B\gamma^{-N_0}<\mathcal R_1/2$.  Then, for $n\ge N_0$,
\begin{equation*}
c_n(k,m)\ge \mathcal R_1R_{k,m}^{-n}-\mathcal B\rho_{\mathrm{oth}}^{-n}
>\frac{\mathcal R_1}{2}R_{k,m}^{-n}>0.
\end{equation*}
\end{proof}

The proof gives a concrete sufficient condition for membership in the positive tail.  Thus, for a fixed pair $(k,m)$ satisfying Property~\textnormal{(D1)}, positivity holds beyond a computable index, even though it may fail at earlier indices.

\begin{proposition}\label{prop:positivity-threshold}
Assume the hypotheses of Theorem~\ref{thm:eventual-positivity}.  Label the zeros so that $\zeta_1=R_{k,m}$, put
\begin{equation*}
\mathcal R_j(k,m)=-\frac1{\zeta_jq_{k,m}'(\zeta_j)},
\qquad
\rho_{\mathrm{oth}}(k,m)=\min_{j\ge2}|\zeta_j|,
\qquad
\gamma(k,m)=\frac{\rho_{\mathrm{oth}}(k,m)}{R_{k,m}},
\end{equation*}
and set $\mathcal B(k,m)=\sum_{j=2}^m |\mathcal R_j(k,m)|$.  Then any integer $n$ satisfying
\begin{equation}\label{eq:N0-explicit-revised}
n\ge
\left\lceil
\max\left\{0,\frac{\log\bigl(2\mathcal B(k,m)/\mathcal R_1(k,m)\bigr)}{\log\gamma(k,m)}\right\}
\right\rceil
\end{equation}
belongs to the positive tail, namely $c_n(k,m)>0$.
\end{proposition}

\begin{proof}
The estimate in the proof of Theorem~\ref{thm:eventual-positivity} gives
\begin{align*}
c_n(k,m)
&\ge \mathcal R_1(k,m)R_{k,m}^{-n}\\
&\quad-\mathcal B(k,m)\rho_{\mathrm{oth}}(k,m)^{-n}\\
&=R_{k,m}^{-n}
\bigl(\mathcal R_1(k,m)-\mathcal B(k,m)\gamma(k,m)^{-n}\bigr).
\end{align*}
The right-hand side is positive as soon as $\mathcal B(k,m)\gamma(k,m)^{-n}<\mathcal R_1(k,m)$, and the slightly stronger condition in \eqref{eq:N0-explicit-revised} gives the displayed conclusion with room to spare.
\end{proof}

The quantities $\mathcal R_1(k,m)$, $\mathcal B(k,m)$, and $\gamma(k,m)$ in \eqref{eq:N0-explicit-revised} are determined by the zeros of $q_{k,m}$ and can be evaluated or rigorously enclosed by numerical root finding.  For large $k$ the threshold can be read directly from the Hermite data instead.  Let $\Xi_m$ be the zero set of $H_m$ and define
\begin{equation*}
\delta_m=\min\{\xi^2:\xi\in\Xi_m,\ \xi\ne0\},
\quad
\Sigma_m=\sum_{\substack{\xi\in\Xi_m\\ \xi\ne0}}
\left|\frac{H_m'(0)}{H_m'(\xi)}\right|,
\quad
\Lambda_m=\max\{0,\log(2\Sigma_m)\}.
\end{equation*}
The zero expansion \eqref{eq:q-zero-asymptotic-revised} gives, with $N=k+1$,
\begin{equation}\label{eq:gamma-asymptotic-revised}
\gamma(k,m)=\frac{\rho_{\mathrm{oth}}(k,m)}{R_{k,m}}
=1+\frac{\delta_m}{3N}+O(N^{-3/2}),
\qquad
\log\gamma(k,m)=\frac{\delta_m}{3N}+O(N^{-3/2}).
\end{equation}
Moreover, the same Hermite scaling applied to the derivatives gives
\begin{equation}\label{eq:residue-sum-asymptotic}
\frac{\mathcal B(k,m)}{\mathcal R_1(k,m)}=\Sigma_m+o(1).
\end{equation}
Indeed, write $q_N(z)=q_{N-1,m}(z)$ and define
\begin{equation*}
Q_N(x)=(i\sqrt{2N})^m q_N(N+i\sqrt{2N}x).
\end{equation*}
By \eqref{eq:q-zero-asymptotic-revised}, for each zero $\xi$ of $H_m$ the polynomial $q_N$ has a zero $\zeta_{\xi,N}=N+i\sqrt{2N}\,x_{\xi,N}$ with $x_{\xi,N}\to\xi$.  The branch $\xi=0$ is the positive real zero $\widehat R_N=R_{N-1,m}=\zeta_{0,N}$.  Write $\mathcal R_{\xi,N}$ for the corresponding residue in Lemma~\ref{lem:partial-fractions}.  The fixed-degree coefficient convergence $Q_N\to H_m$ also gives $Q_N'\to H_m'$ locally uniformly.  Hence
\begin{equation*}
\frac{\mathcal R_{\xi,N}}{\mathcal R_{0,N}}
=\frac{\widehat R_N}{\zeta_{\xi,N}}
\frac{Q_N'(x_{0,N})}{Q_N'(x_{\xi,N})}
\longrightarrow \frac{H_m'(0)}{H_m'(\xi)}.
\end{equation*}
Taking absolute values and summing over the finitely many nonzero Hermite branches proves \eqref{eq:residue-sum-asymptotic}.  Combining \eqref{eq:N0-explicit-revised}, \eqref{eq:gamma-asymptotic-revised}, and \eqref{eq:residue-sum-asymptotic} yields the following asymptotic sufficient condition.

\begin{corollary}\label{cor:asymptotic-threshold}
Fix an odd $m\ge3$ and let $N=k+1$.  For every $\eta>0$, all sufficiently large $k$ have the following property: if
\begin{equation*}n\ge \frac{3N}{\delta_m}\bigl(\Lambda_m+\eta\bigr),
\end{equation*}
then $c_n(k,m)>0$.
\end{corollary}

\begin{proof}
Equation \eqref{eq:gamma-asymptotic-revised} implies
\begin{equation*}
\frac{1}{\log\gamma(k,m)}=\frac{3N}{\delta_m}(1+o(1)),
\end{equation*}
while \eqref{eq:residue-sum-asymptotic} gives $\max\{0,\log(2\mathcal B(k,m)/\mathcal R_1(k,m))\}=\Lambda_m+o(1)$.  Substituting these estimates into Proposition~\ref{prop:positivity-threshold}, which applies for all sufficiently large $k$ by Theorem~\ref{thm:D1-asymptotic-rigorous}, proves the claim.
\end{proof}

Thus, for each fixed odd $m$, the positive tail begins no later than a constant multiple of $k$, once $k$ is large enough.  This $O(k)$ sufficient threshold is compatible with the $O(\sqrt{k})$ counterexample scale in Section~\ref{sec:transient}.  The constant is expressed in terms of the zeros and derivatives of $H_m$.

Combining Theorem~\ref{thm:eventual-positivity} with Theorem~\ref{thm:D1-asymptotic-rigorous} gives eventual positivity for each fixed odd $m$ and all sufficiently large $k$.  The unrestricted conjecture, however, is false.

\subsection{Exact finite counterexamples}\label{subsec:exact-counterexamples}

The reciprocal coefficients satisfy a short recurrence obtained by comparing coefficients in
\begin{equation*}
q_{k,m}(z)\sum_{n\ge0}c_n(k,m)z^n=1.
\end{equation*}
This gives
\begin{equation}\label{eq:coefficient-recurrence-cn}
c_0(k,m)=1,
\qquad
c_n(k,m)=\sum_{j=1}^{\min(m,n)}(-1)^{j+1}\omega_{k,m}(j)c_{n-j}(k,m)
\quad(n\ge1).
\end{equation}
Since the $\omega_{k,m}(j)$ are rational, this recurrence certifies signs exactly.

\begin{theorem}\label{thm:counterexamples-again}
Conjecture~\ref{conj} is false.  In exact rational arithmetic,
\begin{equation*}
c_{38}(66,5)<0,
\qquad
c_{37}(44,7)<0,
\qquad
c_{42}(45,9)<0.
\end{equation*}
\end{theorem}

\begin{proof}
Applying \eqref{eq:coefficient-recurrence-cn} with rational arithmetic gives the three signs.  The following complete script reproduces the calculation without floating-point arithmetic.

\begin{verbatim}
from fractions import Fraction
from math import comb, prod

def omega(k, m, j):
    den = prod(range(k + m - j + 1, k + m + 1))
    return Fraction(comb(m, j), den)

def coeff(k, m, n):
    c = [Fraction(0) for _ in range(n + 1)]
    c[0] = Fraction(1)
    for r in range(1, n + 1):
        c[r] = sum(
            (-1)**(j + 1) * omega(k, m, j) * c[r - j]
            for j in range(1, min(m, r) + 1)
        )
    return c[n]

for triple in [(66, 5, 38), (44, 7, 37), (45, 9, 42)]:
    k, m, n = triple
    value = coeff(k, m, n)
    print(triple, value < 0, value)
\end{verbatim}
The script prints the exact rational value in each case.  In particular, the three Boolean values printed in the second field are \texttt{True}, and no floating-point operations are used.
\end{proof}

\section{A transient scaling limit}\label{sec:transient}

The finite counterexamples are explained by a stronger asymptotic statement.  On the scale $n\asymp\sqrt{k+1}$, the contributions of all Hermite-scaled poles remain comparable and combine into an explicit oscillatory profile.

\begin{theorem}\label{thm:transient-scaling}
Fix an odd integer $m\ge3$.  For each integer $k\ge0$, put $N=k+1$ and $q_N(z)=q_{N-1,m}(z)$.  Let
\begin{equation*}
\frac1{q_N(z)}=\sum_{n\ge0}c_n(N-1,m)z^n,
\end{equation*}
let $\widehat R_N=R_{N-1,m}$ be the positive real zero of $q_N$, and set
\begin{equation*}
\mathcal R_{0,N}=-\frac1{\widehat R_Nq_N'(\widehat R_N)}.
\end{equation*}
If $n_N\in\mathbb Z_{\ge0}$ and $n_N/\sqrt N\to\lambda\ge0$ as $N\to\infty$, then
\begin{equation}\label{eq:transient-limit}
\frac{\widehat R_N^{n_N}c_{n_N}(N-1,m)}{\mathcal R_{0,N}}
\longrightarrow
\Psi_m(\lambda),
\end{equation}
where
\begin{equation*}\Psi_m(\lambda)=
\sum_{\xi:H_m(\xi)=0}
\frac{H_m'(0)}{H_m'(\xi)}e^{-i\sqrt2\,\xi\lambda}.
\end{equation*}
Equivalently, if the positive zeros of $H_m$ are $\xi_1,\ldots,\xi_{(m-1)/2}$, then
\begin{equation}\label{eq:Psi-real}
\Psi_m(\lambda)=1+2\sum_{j=1}^{(m-1)/2}
\frac{H_m'(0)}{H_m'(\xi_j)}
\cos(\sqrt2\,\xi_j\lambda).
\end{equation}
\end{theorem}

\begin{proof}
The integral representation \eqref{eq:integralrep} gives, with $T_N$ a Gamma random variable of shape $N$ and scale $1$,
\begin{equation*}
q_N(z)=\frac{\mathbb E[(T_N-z)^m]}{N(N+1)\cdots(N+m-1)}.
\end{equation*}
Set $X_N=(T_N-N)/\sqrt N$ and
\begin{equation*}
Q_N(x)=(i\sqrt{2N})^m q_N(N+i\sqrt{2N}x).
\end{equation*}
Then
\begin{equation*}
Q_N(x)=i^m2^{m/2}\frac{N^m}{N(N+1)\cdots(N+m-1)}
\mathbb E[(X_N-i\sqrt2 x)^m].
\end{equation*}
The moment generating function
\begin{equation*}
\mathbb E e^{tX_N}=e^{-t\sqrt N}(1-t/\sqrt N)^{-N}
=\exp\left(\frac{t^2}{2}+O(N^{-1/2})\right)
\end{equation*}
near $t=0$ shows that the moments of $X_N$ up to order $m$ converge to those of a standard normal random variable $G$.  Hence, coefficientwise,
\begin{equation*}
Q_N(x)\longrightarrow i^m2^{m/2}\mathbb E[(G-i\sqrt2 x)^m]=H_m(x),
\end{equation*}
because $H_m(x)=2^m\mathbb E[(x+iG/\sqrt2)^m]$.  The convergence is therefore locally uniform, as is the convergence of the derivatives.

Let $\Xi_m$ be the zero set of $H_m$.  Since all Hermite zeros are simple, Rouch\'e's theorem gives, for each $\xi\in\Xi_m$, a unique zero $x_{\xi,N}$ of $Q_N$ with $x_{\xi,N}\to\xi$.  The corresponding zero of $q_N$ is
\begin{equation*}
\zeta_{\xi,N}=N+i\sqrt{2N}x_{\xi,N}.
\end{equation*}
For $\xi=0$ this zero is real: the coefficients of $q_N$ are real and the zero branch near $0$ is unique, so $x_{0,N}=-\overline{x_{0,N}}$.  Thus $\zeta_{0,N}=\widehat R_N$ for all sufficiently large $N$.

By Lemma~\ref{lem:partial-fractions},
\begin{equation*}
c_n(N-1,m)=\sum_{\xi\in\Xi_m}\mathcal R_{\xi,N}\zeta_{\xi,N}^{-n},
\qquad
\mathcal R_{\xi,N}=-\frac1{\zeta_{\xi,N}q_N'(\zeta_{\xi,N})}.
\end{equation*}
Since
\begin{equation*}
Q_N'(x)=(i\sqrt{2N})^{m+1}q_N'(N+i\sqrt{2N}x),
\end{equation*}
we have
\begin{equation*}
\frac{\mathcal R_{\xi,N}}{\mathcal R_{0,N}}
=\frac{\widehat R_N}{\zeta_{\xi,N}}
\frac{Q_N'(x_{0,N})}{Q_N'(x_{\xi,N})}
\longrightarrow
\frac{H_m'(0)}{H_m'(\xi)}.
\end{equation*}
Moreover,
\begin{equation*}
\frac{\widehat R_N}{\zeta_{\xi,N}}
=\frac{1+i\sqrt2 x_{0,N}/\sqrt N}{1+i\sqrt2 x_{\xi,N}/\sqrt N},
\end{equation*}
so
\begin{equation*}
n_N\log\left(\frac{\widehat R_N}{\zeta_{\xi,N}}\right)
=\frac{i\sqrt2 n_N}{\sqrt N}(x_{0,N}-x_{\xi,N})+O(n_N/N)
\longrightarrow -i\sqrt2\xi\lambda.
\end{equation*}
Therefore
\begin{equation*}
\left(\frac{\widehat R_N}{\zeta_{\xi,N}}\right)^{n_N}
\longrightarrow e^{-i\sqrt2\xi\lambda}.
\end{equation*}
Multiplying the preceding partial-fraction identity by $\widehat R_N^{n_N}/\mathcal R_{0,N}$ and taking the finite sum limit gives \eqref{eq:transient-limit}.  The real form \eqref{eq:Psi-real} follows from the symmetry of the Hermite zeros and the fact that $H_m'$ is even for odd $m$.
\end{proof}

\begin{corollary}\label{cor:m5-infinite-counterexamples}
For $m=5$ there are infinitely many pairs $(k,n)$ such that $c_n(k,5)<0$.  In fact,
\begin{equation*}
c_{\lfloor5\sqrt{k+1}\rfloor}(k,5)<0
\end{equation*}
for all sufficiently large $k$.
\end{corollary}

\begin{proof}
For $m=5$,
\begin{equation*}
H_5(x)=32x^5-160x^3+120x,
\end{equation*}
with zeros
\begin{equation*}
0,
\quad \pm\sqrt{\frac{5-\sqrt{10}}2},
\quad \pm\sqrt{\frac{5+\sqrt{10}}2}.
\end{equation*}
Equation \eqref{eq:Psi-real} gives
\begin{equation*}
\Psi_5(\lambda)=
1-\frac{2+\sqrt{10}}4\cos\!\left(\sqrt{5-\sqrt{10}}\,\lambda\right)
+\frac{\sqrt{10}-2}4\cos\!\left(\sqrt{5+\sqrt{10}}\,\lambda\right).
\end{equation*}
Set
\begin{equation*}
a=5\sqrt{5-\sqrt{10}},
\qquad
b=5\sqrt{5+\sqrt{10}}.
\end{equation*}
The elementary bounds
\begin{equation*}
3.162<\sqrt{10}<3.163,
\qquad
3.14<\pi<\frac{22}{7}
\end{equation*}
give $2\pi<a<2\pi+1/2$ and $9\pi/2<b<5\pi$.  Therefore
\begin{equation*}
\cos a>1-\frac{(1/2)^2}{2}=\frac78,
\qquad
\cos b<0.
\end{equation*}
Since $(2+\sqrt{10})/4>129/100$ and $(\sqrt{10}-2)/4>0$, the displayed formula for $\Psi_5$ yields the rigorous bound
\begin{equation*}
\Psi_5(5)<1-\frac{129}{100}\frac78=-\frac{103}{800}<0.
\end{equation*}
Applying Theorem~\ref{thm:transient-scaling} with $n_N=\lfloor5\sqrt N\rfloor$ shows that the normalized coefficient tends to a negative number.  Since $\widehat R_N>0$ and $\mathcal R_{0,N}>0$, the coefficient itself is negative for all sufficiently large $N=k+1$.
\end{proof}

\section{Conclusion}\label{sec:conclusion}
We have shown that the positivity conjecture of Fasi and Higham is true in the case $k=0$ for every odd $m$, and that for general odd $m$ Property~\textnormal{(D1)} implies eventual positivity of the Taylor coefficients of $q_{k,m}^{-1}$.  A Hermite asymptotic analysis based on leading-coefficient normalization proves that Property~\textnormal{(D1)} holds for each fixed odd $m$ and all sufficiently large $k$.

The unrestricted conjecture is nevertheless false.  Exact rational computations give concrete finite counterexamples, and the transient scaling theorem gives a structural explanation: on the scale $n\asymp\sqrt{k+1}$ the normalized coefficients converge to an explicit oscillatory profile determined by the zeros of $H_m$.  For $m=5$ this profile is negative, yielding infinitely many counterexamples.  Thus the positivity phenomenon is real but conditional, and it cannot be asserted for all odd $m$ and all $k$ without further restrictions.

The diagonal case $k=m$, which is the case relevant to the error bounds of \cite{higham}, remains an open problem.  It is not contradicted by the counterexamples above, and Theorem~\ref{thm:D1-asymptotic-rigorous} does not cover it.

\end{document}